\documentclass[11pt]{article}
\usepackage{amsmath}
\usepackage[]{graphicx}

\newcommand{\bmath}{\begin{displaymath}}
\newcommand{\emath}{\end{displaymath}}
\newcommand{\beq}{\begin{equation}}
\newcommand{\eeq}{\end{equation}}
\newcommand{\beqa}{\begin{eqnarray}}
\newcommand{\eeqa}{\end{eqnarray}}

\begin{document}

\title{On ``jamitons,'' self-sustained nonlinear traffic waves}
\author{M.~R.~Flynn$^1$, A.~R.~Kasimov$^2$, J.-C.~Nave$^2$, R.~R.~Rosales$^2$, B.~Seibold$^2$ \\ \\ 
$^1$ Dept. of Mechanical Engineering, Univ. of Alberta, \\ Edmonton, AB, T6G 2G8 Canada \\ \\
$^2$ Dept. of Mathematics, Massachusetts Inst. of Technology, \\ 77 Massachusetts Avenue, Cambridge, MA 02139, USA} 

\date{\today}
\maketitle

\begin{abstract}
``Phantom jams,'' traffic blockages that arise without apparent cause, have long frustrated transportation scientists. Herein, we draw a novel homology between phantom jams and a related class of self-sustained transonic waves, namely detonations. Through this analogy, we describe the jam structure; favorable agreement with reported measurements from congested highways is observed. Complementary numerical simulations offer insights into the jams' development. Our results identify conditions likely to result in a dangerous concentration of vehicles and thereby lend guidance in traffic control and roadway design.
\end{abstract}


\vspace{1cm}

Continuum traffic models whereby individual vehicles are modeled collectively have enjoyed a rich history since their popularization by M.~J. Lighthill and G.~B. Whitham. These models are able to explain dynamics both on open roadways and in the vicinity of traffic lights and also elucidate ``phantom jams,'' blockages that arise without obstructions \cite{helbing2001}.
Notwithstanding this venerated record, the possibility of self-sustained disturbances consisting of a shock matched to a transonic flow has received inadequate attention despite recent experimental evidence that such ``jamitons'' are ubiquitous \cite{sugiyama_fukui_kikuchi_hasebe_nakayama_nishinari_tadaki_yukawa2008}.
We demonstrate herein that jamitons are structurally similar to detonation waves \cite{fickett_davis1979}. 
This insight, fascinating in its own right, lends helpful guidance in roadway design. 

Continuum traffic models typically consist of coupled pair of partial differential equations describing the mass and momentum of the flow \cite{helbing2001}. Using subscripts to indicate differentiation in time and space, this coupled pair is expressed, generically, as
\beq
\rho_t+(\rho u)_x=0\, ,
\label{eq:mass}
\eeq
\beq
u_t+u u_x+\frac{p_x}{\rho}=\frac{1}{\tau}\left(\tilde{u}-u\right)\, ,
\label{eq:mom}
\eeq
where $\tau$ is a relaxation timescale and $u$ and $\rho$ are, respectively, the local traffic speed and density. From the discussion of Aw \& Rascle \cite{aw_rascle2000}, 
the traffic pressure, $p$, increases with density. The desired speed, $\tilde{u}$, decreases with density and vanishes as $\rho$ approaches the maximal density, $\rho_M$.  

\begin{figure}
\includegraphics[width=.95\textwidth]{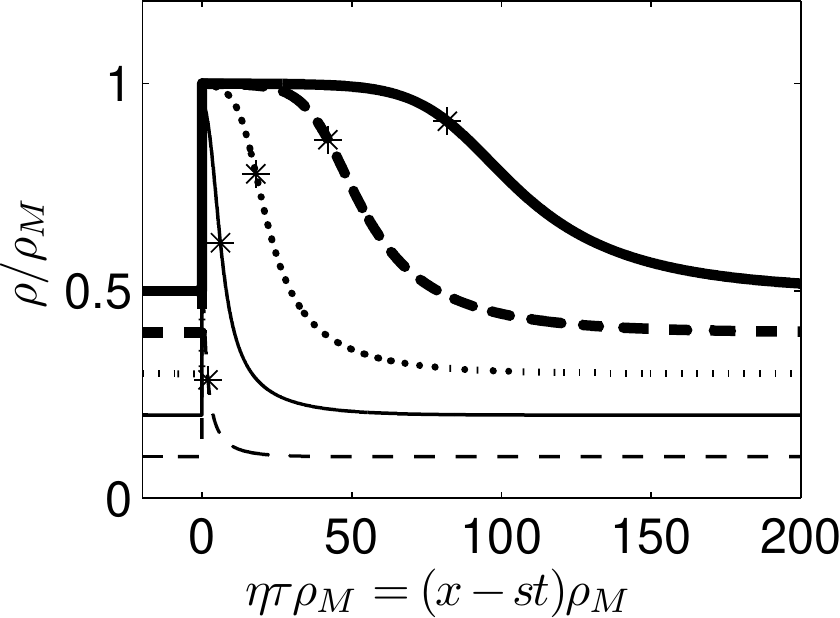}
\caption{Profiles of the non-dimensional traffic density for
$(\rho_-/\rho_M,\,u_-/\tilde{u}_0)=(0.1,\,0.9),\,(0.2,\,0.8),\,\ldots\,(0.5,\,0.5)$,
where $\rho_-$ and $u_-$ denote the traffic density and speed far up- or downstream of the shock.
Here $p_x/\rho=\beta\rho_x/(\rho_M-\rho)$ and $\tilde{u}=\tilde{u}_0(1-\rho/\rho_M)$ with
$\beta=10$\,m$^2$/s$^2$, $\rho_M=0.2$\,m$^{-1}$, $\tilde{u}_0=20$\,m/s and $\tau=5$\,s.
Stars denote the sonic point.}
\label{fg:profiles}
\end{figure}

In the experiment of Sugiyama et al.~\cite{sugiyama_fukui_kikuchi_hasebe_nakayama_nishinari_tadaki_yukawa2008}, 22 vehicles were placed on closed 0.23\,km long track; the distance between adjacent vehicles was initially constant. Instructing drivers to travel at 30\,km/h, researchers soon discerned spatial inhomogeneities in $\rho$ and $u$. Most notably, a propagating density bulge, approximately six vehicles thick, was observed to travel against the direction of traffic flow. Away from this jam, vehicles moved unencumbered. Each vehicle decelerated as it entered the jam, was momentarily delayed, then accelerated out the other side. 

Motivated by these observations, self-sustained solutions are sought (on an infinite road) in terms of the self-similar variable $\eta=(x-st)/\tau$ where $s\,(<u)$ is the traveling wave speed. In terms of $\eta$, (\ref{eq:mass}) and (\ref{eq:mom}) may be combined to yield
\beq
\frac{\mathrm{d}u}{\mathrm{d}\eta}=\frac{(u-s)(\tilde{u}-u)}{(u-s)^2-c^2}\, .
\label{eq:deton}
\eeq

\begin{figure}
\includegraphics[width=.99\textwidth]{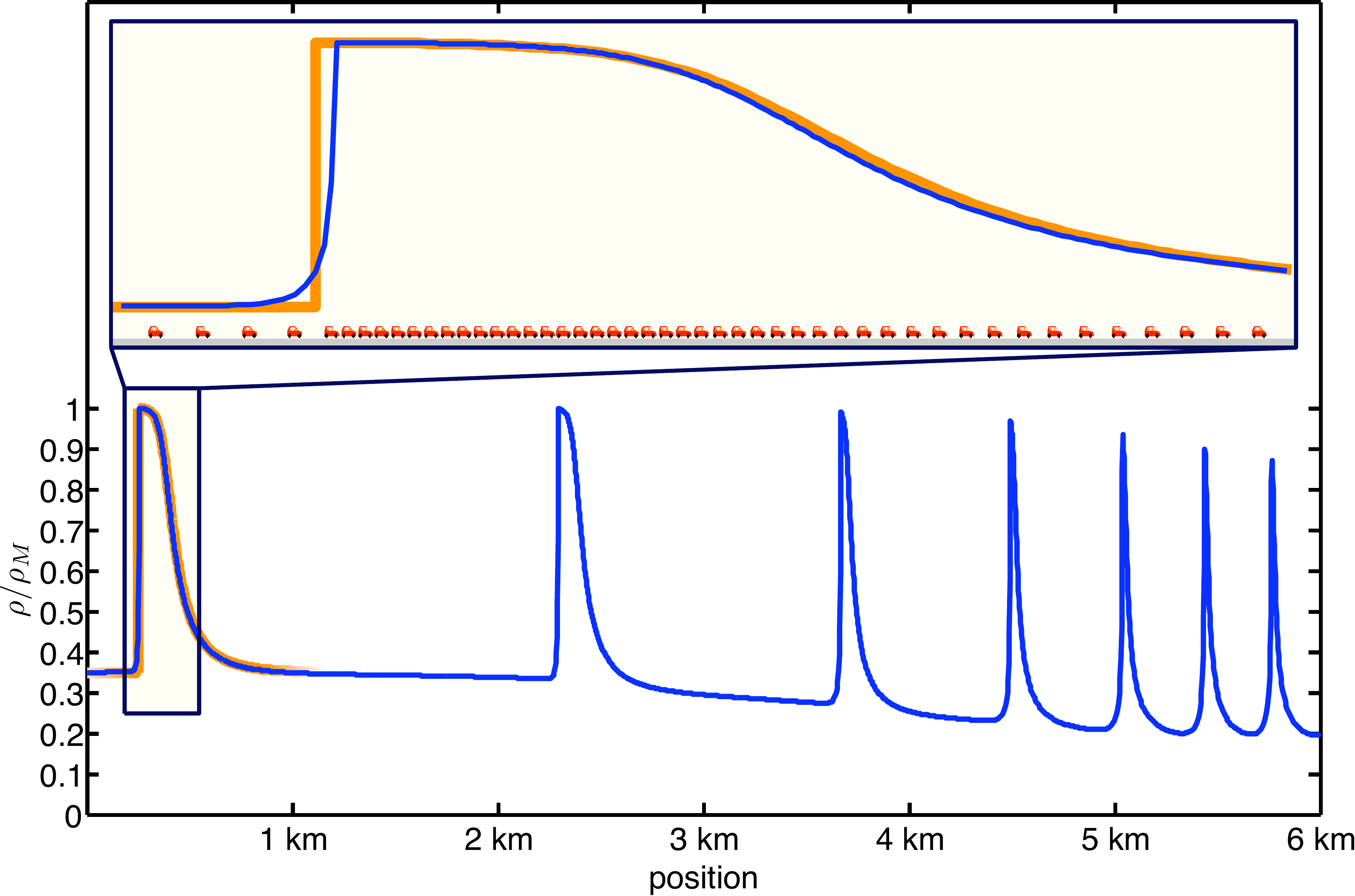}
\caption{Multiple jamitons. Orange and blue curves show, respectively, the theoretical and
numerical solutions. The difference is visible only in the expanded inset, which also presents
a representative distribution of individual vehicles along the bottom.}
\label{fg:simulation}
\end{figure}

Here $c \equiv (\mathrm{d}p/\mathrm{d}\rho)^{1/2}$ is analogous to the sound speed from compressible flow. The solution of (\ref{eq:deton}) follows from analyses in disciplines as distinct as combustion \cite{fickett_davis1979}, 
astrophysics \cite{chakrabarti1990}, 
and arterial blood flow \cite{jensen_pedley1989}. 
A sonic point is an event horizon where, in the frame of the traveling wave, the vehicle speed matches the speed of small disturbances, and is defined by $u-s=c$. The speed and structure of the jamiton is obtained by integration of (\ref{eq:deton}) under the condition that $\mathrm{d}u/\mathrm{d}\eta$ remain finite at the sonic point where both the numerator and denominator of (\ref{eq:deton}) simultaneously vanish. The jamiton is a self-sustained wave because the region between the shock and the sonic point is decoupled from the flow downstream of the sonic point in direct analogy with self-sustained (Chapman-Jouguet) detonation waves \cite{fickett_davis1979}.

Solutions are presented in figure~\ref{fg:profiles}, which shows profiles of the traffic density. Consistent with Sugiyama et al.~\cite{sugiyama_fukui_kikuchi_hasebe_nakayama_nishinari_tadaki_yukawa2008}, jamitons connect up- and downstream states of equal density and speed and arise only when, as in
figure~\ref{fg:profiles}, the background density exceeds a critical threshold. There is a striking correspondence between the curves of
figure~\ref{fg:profiles} and the detonation wave depicted in figure~2.1\,a of Fickett \& Davis \cite{fickett_davis1979}.

As with the roll waves observed in gutters on rainy days, numerical solution\footnote[1]{Details of the numerical method are briefly reviewed in Appendix A. Further detail will be presented in a forthcoming publication. Additional examples of numerical output are given
at \texttt{http://math.mit.edu/projects/traffic}.} of (\ref{eq:mass}) and (\ref{eq:mom}) shows the development of a train of jamitons.
Figure~\ref{fg:simulation} indicates that individual waves from this train are well-described by the foregoing theoretical analysis even when, as with the left-most jamiton, the waves have not completely saturated. 

\begin{figure}
\begin{minipage}[t]{.45\textwidth}
\centering
\includegraphics[width=\textwidth]{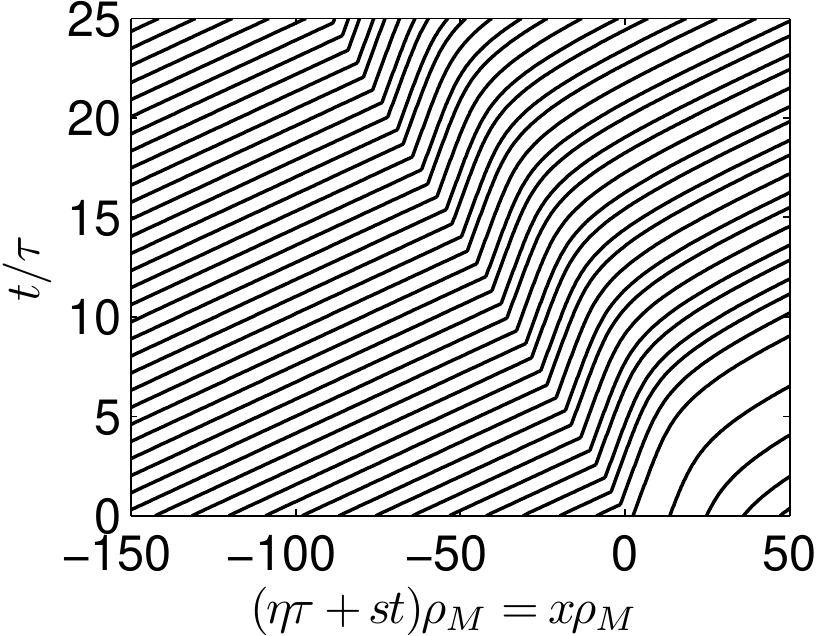}
\caption{Vehicle trajectories for $(\rho_-/\rho_M,\,u_-/\tilde{u}_0)=(0.35,\,0.65)$.}
\label{fg:trajectories_theory}
\end{minipage}
\hfill
\begin{minipage}[t]{.50\textwidth}
\centering
\includegraphics[width=\textwidth]{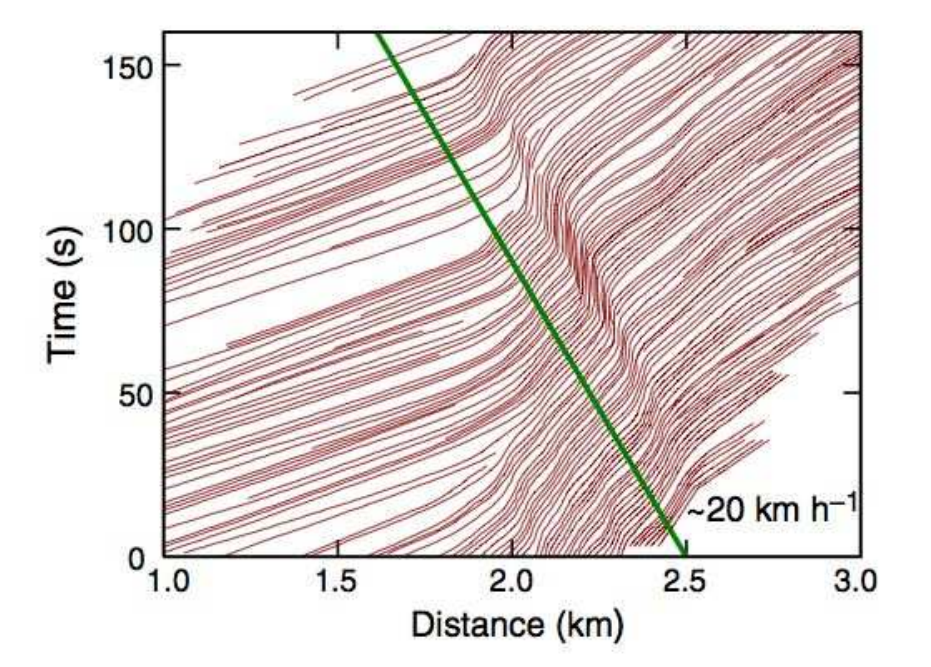}
\caption{Measured vehicle trajectories due to J.~Treiterer. The green line shows the experimental
measurement of \cite{sugiyama_fukui_kikuchi_hasebe_nakayama_nishinari_tadaki_yukawa2008}
(Figure reproduced from \cite{sugiyama_fukui_kikuchi_hasebe_nakayama_nishinari_tadaki_yukawa2008}
with permission).}
\label{fg:trajectories_experiment}
\end{minipage}
\end{figure}

Vehicle trajectories are shown in figure~\ref{fg:trajectories_theory}. Left-to-right moving vehicles travel at constant speed $u_-$ until they encounter the upstream-propagating jam. Thereafter, a rapid deceleration, followed by a gradual acceleration, is observed.
Data from congested highways show similar trends (figure~\ref{fg:trajectories_experiment}). 

Descriptions of jamitons provide helpful design guidance as they identify conditions likely to produce dangerous vehicle concentrations. Such situations may be avoided by judicious selection of speed limits, carrying capacities, etc.

\vspace*{0.5cm}
\noindent
{\it Acknowledgments}\\
Funding generously provided by NSF (DMS-0813648) and the AFOSR Young Investigator Program (FA9550-08-1-0035).

\appendix
\section{Numerical method}

To investigate the time evolution of density and speed anomalies within a stream of vehicles, the governing equations (1) and (2) are solved numerically using a mesh-free, Lagrangian particle method \cite{monaghan1988} 
calibrated to accommodate rapid changes in $\rho$ and $u$ in the vicinity of (numerically-smeared) shocks. The continuum of vehicles is represented as a set of $\mathcal{O}(10^3)-\mathcal{O}(10^4)$ discrete particles whose motions are governed by the combination of forces given on the right-hand-side of (2). Individual particles neither appear or disappear; consequently the mass continuity equation (1) is satisfied automatically. The number of particles exceeds by two or more orders of magnitude the number of vehicles. Therefore, although there is a natural connection between the particle and vehicle density, the numerical method is fundamentally ``macroscopic'' in that the reaction of individual drivers is not explicitly modeled. 

\bibliographystyle{unsrt}
\bibliography{fluid_refs}

\end{document}